# The domain of the Fourier integral

V. N. Tibabishev

**We consider the problem of determining the Fourier integral in the Hilbert space of square integrable functions. Fourier integral is the scalar product of two functions belonging to the Hilbert space of square integrable functions and the Hilbert space of almost periodic functions. Scalar product for different Hilbert spaces defined at the intersection of these spaces, which contains only one zero element. Therefore, the Fourier integral is not defined in the Hilbert space of square integrable functions with nonzero norm.**

Keywords: Hilbert spaces of almost periodic functions, Fourier transform

**Introduction**

In automatic control theory has found wide application of the Fourier integral. Known [1, s.397, 412] that the Fourier integral is defined in the space of absolutely integrable functions on the whole line $L_1(-\infty,+\infty)$ and in the Hilbert space of square integrable functions on the whole line $L_2(-\infty,+\infty)$. The functions performed by the Fourier integral, is interpreted as an infinite sum of harmonic oscillations with infinitely close frequencies and with an infinitely small amplitudes [2, p.16]. In addition, in a Hilbert space $L_2(-\infty,+\infty)$ there exists an orthonormal system of functions, such as Hermite [1, pp. 375]. In this case, the same element of the Hilbert space $L_2(-\infty,+\infty)$ can be represented not only by the Fourier integral, but also a Fourier series in the orthonormal system, for example, the Hermite functions. Dual representation of one element belonging to the same Hilbert space, leads to contradictions. In the representation of the Fourier space element is interpreted by a countably infinite sum of projections of a countable orthonormal basis $L_2(-\infty,+\infty)$. In the representation of the Fourier integral, the same element in the same space is interpreted by an uncountable amount of harmonic vibrations, which do not belong to a Hilbert space $L_2(-\infty,+\infty)$. It is known [3, pp. 31] that what would have been the power of an orthonormal system in a Hilbert space $H$, every vector $f \in H$ has at most a countable set of nonzero projections on the elements of an orthonormal system. In this regard, seeks to determine the causes of the contradictions in the theory of Fourier integral, as defined in Hilbert space $L_2(-\infty,+\infty)$.

## 2. Two models of representation of non-periodic functions

Fourier integral is introduced as a generalization of the Fourier series, defined for periodic functions in the event of non-periodic functions [1, pp. 393]. It is known [2, pp. 12] that the function is called a periodic function, if any, the equality

$$f(t) = f(t + nT), \qquad (1)$$

where $T$ - constant, $n$ - any integer positive or negative. There are two ideas of generalization of periodic functions represented by Fourier series of harmonic functions with multiple frequencies, the class of aperiodic functions. The first idea of generalization of periodic functions to the class of aperiodic functions is assumed that the infinite increase in the period of the fundamental harmonic $T \to \infty$ will be obtained non-periodic function [2, pp. 15]. This idea contradicts the definition of a periodic function (1). Known [4, pp. 66] that an infinitely large value is not constant, but variable, and therefore infinitely large value of the period cannot be substituted in formula (1) by definition. If you still do this substitution, we find that for any real value of the argument $t$ we have $f(t) = f(t + n\infty) = f(\infty)$. This equation holds for a variety of functions $f_i(t) = C_i$, where $C_i$ - arbitrary constants, $i = 1,2,3,...$. If $C_i \neq 0$, then such a set of functions does not belong to $L_1(-\infty,+\infty)$ nor to $L_2(-\infty,+\infty)$. Only function $f(t) \equiv 0$ is an element of space $L_1(-\infty,+\infty)$ and space $L_2(-\infty,+\infty)$. This leads to the conclusion that the first idea of generalized periodic functions on the class of aperiodic functions will lead to transformation, defined in the space $L_1(-\infty,+\infty)$ and $L_2(-\infty,+\infty)$ only for one function $f(t) \equiv 0$.

Another idea is a generalization of periodic functions to the class of non-periodic functions is that the periodic function represented by Fourier series with multiple frequencies, becomes non-periodic function if multiple frequencies replaced by arbitrary incommensurate frequencies. For example, the sum of two harmonic oscillations with finite period, but with differing frequencies, for example, $\omega_1$ and $\sqrt{2}\omega_1$ leads to a clearly aperiodic fluctuations [2, pp. 18]. Sum of harmonic vibrations with arbitrary frequencies studied in the theory of almost periodic functions. There are various definitions and generalizations of almost periodic functions. Among them Besicovitch almost periodic functions [5, p. 222], for which the usual quantity

$$D_{B^p}[f,g] = \left\{M\left\{|f-g|^p\right\}\right\}^{1/p}, \text{ где } M\{\cdot\} = \lim \frac{1}{T}\int_{-T/2}^{+T/2}\{\cdot\}dt \text{ при } T\to\infty \quad (2)$$

is the $B^p$ distance. Space of functions summable with a $p$ second degree in every finite interval with such a certain distance, called $B^p$ - space.

The function $f(t)$ is called almost periodic in Besicovitch $p$ order if there is a sequence of finite trigonometric sums $S_1(t), S(t)_2,....S_n(t),...$ for which

$$\lim D_{B^p}[f(t), S_n(t)] = 0 \text{ при } n\to\infty. \quad (3)$$

Spaces $B^p$ ($p \geq 1$) are full. If $p = 2$ space $B^2$ is a Hilbert space in which every almost periodic function $f(t) \in B^2$ can be associated with a Fourier series

$$f(t) \sim \sum_n A_n \exp(j\omega_n t), \text{ где } A_n = M\{f(t)\exp(-j\omega_n t)\}. \quad (4)$$

In the real world can not get the motion of material bodies with nonzero rest mass with infinite acceleration. Therefore, the functions $f(t)$ that describe the real dynamic processes that can only have continuous first derivatives $\dot{f}(t)$ on the whole line. It is known [6, pp. 461] that a function with continuous first derivatives are convergent Fourier series. From the condition (3) that such functions are almost periodic functions in the sense of Besicovitch and belong to the subspace $W_2(-\infty,+\infty) \subset B^2(-\infty,+\infty)$ for which the series (4) is uniformly convergent.

Thus, there are two Hilbert spaces to describe two models of non-periodic functions. Non-recurrent functions on the first model will be denoted by $x(t)$, which presented a Fourier integral in a Hilbert space $L_2(-\infty,+\infty)$. Nonperiodic functions of the second model, denoted by $y(t)$, which seem convergent Fourier series (4) with differing frequencies in a Hilbert space $W_2(-\infty,+\infty)$.

### 3. The domain of the Fourier integral in the space $L_2(-\infty,\infty)$

In the proof of the theorem Plancherel original functional is the functional defining the inner product in Hilbert space $L_2$ with finite measure [1, pp. 412]. Scalar product $(x,y)_H$ defined in a Hilbert space $H$, is only for items $x, y \in H$. If the elements $x$ and $y$ belong to different Hilbert spaces, for example, $x \in L_2$ and $y \in W_2$ then the inner product $(x,y)_{L_2}$ can only be determined on a nonempty subset $Q = L_2 \bigcap W_2 \subset L_2$. In this case, the scalar product $(x, y)_{L_2}$ of elements $x$ and $y$ is not defined on the entire space $L_2$, but only for elements belonging to the subspace $Q \subset L_2$.

In the functional Plancherel inner product is for the two functions $x \in L_2(-\infty,+\infty)$ and $y(t) = \exp(-j\omega t)$, moreover $y \in W_2(-\infty,+\infty)$. Let finite functions, $x_T$ and $y_T$ coincide with the functions $x$ and $y$ on a finite interval $-T < t < +T$ and zero outside this interval. In this case $x_T \in L_2(-T,+T)$ as well $y_T \in W_2(-T,+T)$. Define the subspace of finite measure $Q(-T,+T) = L_2(-T,+T) \bigcap W_2(-T,+T)$. Since nondecreasing integrable $y(t)$ at infinity in finite limits of integration as space $L_2(-T,+T)$, then

$$Q(-T,+T) = L_2(-T,+T) \bigcup W_2(-T,+T). \qquad (5)$$

In this case, the metric $L_2(-T,+T)$ is well-defined inner product for functions of compact support

$$(x_T, \bar{y}_T)_Q = \int_{-T}^{+T} x_T(t) \exp(-j\omega t) dt.$$

In the limit $T \to \infty$ of all functions $y \in W_2(-\infty,+\infty)$ with non-zero norm in the metric $W_2(-\infty,+\infty)$ does not belong to a Hilbert space $L_2(-\infty,+\infty)$. Therefore $T \to \infty$, when condition (5) is not satisfied, but a subset $Q(-\infty,+\infty)$ is determined by the suppression of spaces $Q(-\infty,+\infty) = L_2(-\infty,+\infty) \bigcap W_2(-\infty,+\infty)$. This subspace contains only one zero element. It follows that the Fourier integral is not defined for all elements $x \in L_2(-\infty,+\infty)$, as stated in the Plancherel theorem, and only one zero element $x(t) \equiv 0$.

### 4. Conclusion

Contradictions that occur when using the Fourier integral for the elements belonging to a Hilbert space $L_2(-\infty,+\infty)$, due to the incorrect assertion that the Fourier integral is applicable to all elements of this space. Fourier transform is defined for a set of functions that form the intersection of the Hilbert space $L_2(-\infty,+\infty)$ and another Hilbert space with an uncountable harmonic basis $W_2(-\infty,+\infty)$. The intersection of these spaces $L_2(-\infty,+\infty) \bigcap W_2(-\infty,+\infty)$ has only one zero element. In the Hilbert space with infinite measure $L_2(-\infty,+\infty)$ each element of this space with a nonzero norm is a generalized Fourier series in the orthonormal basis, for example, the normalized Hermite functions. Representation of generalized Fourier series is unique, since for functions with non-zero norm of the Fourier integral is not defined. Fourier integral in a Hilbert space $L_2(-\infty,+\infty)$ is only defined for the zero element. Dual representation of the zero element side and the Fourier integral is correct, since the zero element at the same time belongs to the Hilbert space of functions decreasing at infinity $L_2(-\infty,+\infty)$ and the Hilbert space with infinite measure of non-decreasing function at infinity $W_2(-\infty,+\infty)$.

Functions that describe the real dynamic processes have continuous first derivatives and are represented by convergent Fourier series in harmonic functions. Series with differing frequencies of harmonic components describing aperiodic processes. Nonperiodic processes have discrete spectra. The notion of a continuous spectrum can be introduced formally only in the trivial case when there is a function of the Fourier integral is zero on the entire line.